\documentclass{article}

\usepackage{amssymb,amsmath,latexsym,amsthm}
\usepackage[english]{babel}

\usepackage{hypbmsec}

\newtheorem{theorem}{Theorem}
\newtheorem{lemma}{Lemma}

\theoremstyle{definition}

\begin{document}

\title{On a Lower Bound for $\|(4/3)^k\|$}

\author{Yury PUPYREV}


\maketitle


\begin{abstract}
We prove, that
$$
\biggl\|\biggl(\frac43\biggr)^k\biggr\|
>\biggl(\frac49\biggr)^k\qquad\text{for}\quad k\ge6,
$$
where $\|\cdot\|$ is a distance to the nearest prime.
\end{abstract}

\section{Introduction}
\label{ss1}

In 1994 Bennett \cite{ref1} considered
a generalization of Waring's problem,
namely, a problem on the order $g_N(k)$
of the additive basis
$$
S_N^{(k)}=\{1^k,N^k,(N+1)^k,\dots\},\qquad N\ge2,
$$
of the set of positive integers.
He established the following estimates for\linebreak $\|(1+1/N)^k\|$:
$$
\biggl\|\biggl(1+\frac1{N}\biggr)^k\biggr\|
>3^{-k}\qquad\text{for}\quad4\le N\le k\cdot3^k
$$
and with their help obtained the representation
$$
g_N(k)=N^k+\biggl[\biggl(1+\frac1{N}\biggr)^k\biggr]-2
$$
for $4\le N\le(k+1)^{(k-1)/k}-1$.

He concluded that he needed the inequality
\begin{equation}
\label{eqI}
\biggl\|\biggl(\frac43\biggr)^k\biggr\|
>\biggl(\frac49\biggr)^k\qquad\text{for}\quad k\ge6,
\end{equation}
for the representation
$$
g_3(k) = 3^k + \biggl[\biggl(\frac43\biggr)^k\biggr] - 2.
$$

In 2007  Zudilin~\cite{ref6},
by modifying Baker's construction, namely,
by considering Pad\'e approximations to
the remainder of the series
$$
\frac1{(1-z)^{m+1}}=\sum_{n=0}^\infty\binom{m+n}{m}z^n
$$
and by receiving sharp estimates for the $p$-adic orders
of the arising binomial coefficients,
arrived at the bound
$$
\biggl\|\biggl(\frac43\biggr)^k\biggr\|
>0.4914^k\qquad\text{for}\quad k\ge K,
$$
where $K$ is an effective constant.

In 2009 this author \cite{ref4} received
an exact value of $K$, but it was
too big for checking \eqref{eqI}
for $6\le k< K$.

In this paper using the same method
as Zudilin, but with another set of parameters,
we receive the bound \eqref{eqI}
for $k\ge17\,545\,718$, and check
it using a lemma similar to \cite[Proposition~1]{ref2},
with software for remaining $k$.

Thus, we prove the following result.

\medskip

\begin{theorem}
\label{t1}
We have
$$
\biggl\|\biggl(\frac43\biggr)^k\biggr\|
>\biggl(\frac49\biggr)^k\qquad\text{holds for}\quad k\ge6.
$$
\end{theorem}

\section{Pad\'e Approximations}
\label{ss2}
Following \cite{ref6}, we fix
two positive integers $a$ and $b$,
with $3a\le b$, and write
\begin{align}
\nonumber
\biggl(\frac43\biggr)^{2(b+1)}
&=\biggl(\frac{16}9\biggr)^{(b+1)}
=2^{b+1}\biggl(1+\frac18\biggr)^{-(b+1)}
\\
\nonumber
&=2^{b+1}\sum_{l=0}^\infty\binom{b+l}b
\biggl(\frac18\biggr)^l(-1)^l
\\
\nonumber
&=2^{b+1}2^{-3a}\sum_{l=0}^\infty\binom{b+l}b
2^{3(a-l)}(-1)^{-l}
\\
\nonumber
&=2^{b-3a+1}(-1)^a\sum_{l=0}^\infty\binom{b+l}b
2^{3(a-l)}(-1)^{(a-l)}
\\
\nonumber
&\equiv2^{b-3a+1}(-1)^a\sum_{l=a}^\infty\binom{b+l}b
2^{3(a-l)}(-1)^{(a-l)}\quad(\operatorname{mod}\mathbb Z)
\\
&\equiv2^{b-3a+1}(-1)^a\sum_{\nu=0}^\infty\binom{a+b+\nu}b
\biggl(-\frac18\biggr)^\nu\quad(\operatorname{mod}\mathbb Z).
\label{eq6}
\end{align}
So, we are going to consider Pad\'e approximations
to the function
\begin{equation}
\label{eq7}
F(z)=F(a,b;z)=\sum_{\nu=0}^\infty\binom{a+b+\nu}bz^\nu.
\end{equation}
For any positive integer
$n\le b$ we find \cite{ref6}
\begin{align}
\label{eq8q}
Q_n(z^{-1})&=\frac{(a+b+n)!}
{(a+n-1)!\mspace{2mu}n!\mspace{2mu}(b-n)!}
\int_0^1t^{a+n-1}(1-t)^{b-n}(1-z^{-1}t)^n\,dt
\\
\intertext{and}
\label{eq8p}
R_n(z)&=\frac{(a+b+n)!}
{(a+n-1)!\mspace{2mu}n!\mspace{2mu}(b-n)!}
\cdot z^n
\\
&\qquad\qquad
\times\int_0^1t^n(1-t)^{a+n-1}(1-zt)^{-(a+b+n)-1}\,dt
\nonumber
\end{align}
such that
\begin{equation}
\label{eq9}
Q_n(z^{-1})F(z)=P_n(z^{-1})+R_n(z),
\end{equation}
is performed with
polynomial $P_n(x)\in\mathbb Z[x]$,
$\deg P_n\le n-1$.

\section{Arithmetic argument}
\label{ss3}
For every prime $p>\sqrt{a+b+n}$
we set
\begin{align*}
e_{p,n}&=\min_{\mu\in\mathbb Z}\biggl(
-\biggl\{-\frac{a+n}p\biggr\}
+\biggl\{-\frac{a+n+\mu}p\biggr\}
+\biggl\{\frac{\mu}p\biggr\}
\\
&\qquad\hphantom{=\min_{\mu\in\mathbb Z}\biggl(}
-\biggl\{\frac{a+b+n}p\biggr\}
+\biggl\{\frac{a+b+\mu}p\biggr\}
+\biggl\{\frac{n-\mu}p\biggr\}\biggr),
\\
e'_{p,n}&=\min_{\mu\in\mathbb Z}\biggl(
-\biggl\{\frac{a+n+\mu}p\biggr\}
+\biggl\{\frac{a+n}p\biggr\}
+\biggl\{\frac{\mu}p\biggr\}
\\
&\qquad\hphantom{=\min_{\mu\in\mathbb Z}\biggl(}
-\biggl\{\frac{a+b+n}p\biggr\}
+\biggl\{\frac{a+b+\mu}p\biggr\}
+\biggl\{\frac{n-\mu}p\biggr\}\biggr),
\end{align*}
and for
$$
\Phi=\Phi(a,b,n)=\prod_{p>\sqrt{a+b+n}}p^{e_{p,n}},\qquad
\Phi'=\Phi'(a,b,n)=\prod_{p>\sqrt{a+b+n}}p^{e'_{p,n}},
$$
by lemmas 3 and 4 in \cite{ref6}, we have
\begin{align*}
\Phi^{-1}Q_n(x),\,
\Phi^{-1}P_n(x)&\in\mathbb Z[x],
\\
\Phi'^{-1}(n+1)Q_{n+1}(x),\,
\Phi'^{-1}(n+1)P_{n+1}(x)&\in\mathbb Z[x].
\end{align*}

\section({A Bound for \$||(4/3)k||\$}){A Bound for $\|(4/3)^k\|$}
\label{ss4}
For $a$, $b$, and $n$ we write
$$
a=\alpha m,\quad b=\beta m,\quad n=\gamma m,\qquad
m\in\mathbb N.
$$
Our aim is to find a lower bound for
the absolute value of $\varepsilon_k$,
where
$$
\biggl(\frac43\biggr)^k=M_k+\varepsilon_k,\qquad
M_k\in\mathbb Z,\quad
0<|\varepsilon_k|<\frac12\mspace{2mu}.
$$
For $k\ge3$ we write $k=2(\beta m+1)+j$
with positive integers $m$ and $j<2\beta$.
We multiply both sides of \eqref{eq9} by
$\widetilde\Phi^{-1}2^{b-3a+1+2j}(-1)^a$
(where
$\widetilde\Phi$ is equal to $\Phi$ or
to $\Phi'/(n+1)$; we discuss this
choice in what follows) and put $z=-1/8$:
\begin{equation}
\label{eq14}
\begin{split}
&Q_n(-8)\widetilde\Phi^{-1}3^j
\cdot\biggl(\frac43\biggr)^j2^{b-3a+1}(-1)^a
F\biggl(a,b,-\frac18\biggr)
\\
&\qquad
=P_n(-8)\widetilde\Phi^{-1}2^{b-3a+1+2j}(-1)^a
+R_n\biggl(-\frac18\biggr)
\widetilde\Phi^{-1}2^{b-3a+1+2j}(-1)^a.
\end{split}
\end{equation}
From \eqref{eq6} and \eqref{eq7} one can find that
$$
\biggl(\frac43\biggr)^j2^{b-3a+1}(-1)^a
F\biggl(a,b,-\frac18\biggr)
\equiv\biggl(\frac43\biggr)^{2(b+1)+j}\quad
(\operatorname{mod}\mathbb Z)
=\biggl(\frac43\biggr)^k,
$$
so the left-hand side can be written as
$M'_k+\varepsilon_k$ and one can
rewrite~\eqref{eq14} as
\begin{equation}
\label{eq15}
Q_n(-8)\widetilde\Phi^{-1}3^j\cdot\varepsilon_k
=M''_k+R_n\biggl(-\frac18\biggr)
\widetilde\Phi^{-1}2^{b-3a+1+2j}(-1)^a.
\end{equation}

At this point we should check
if the number $M''_k$ is distinct from zero.
Lemma 2 in \cite{ref6} guarantees
that for $n$ or for $n+1$ we
have $M''_k\ne0$. So
$\widetilde\Phi=\Phi$ if \eqref{eq15} holds for $n$ with
$M''_k\ne0$, and $\widetilde\Phi=\Phi'/(n+1)$ if
\eqref{eq15} holds for $n+1$ with
$M''_k\ne0$. (This odd way of working with $\widetilde\Phi$
becomes more understandable once we have determined
$a$, $b$, and $n$.)

So, assuming that
\begin{equation}
\label{eq100}
\biggl|R_n\biggl(-\frac18\biggr)
\widetilde\Phi^{-1}2^{b-3a+1+2j}(-1)^a\biggr|
<\frac23\mspace{2mu},
\end{equation}
from \eqref{eq15} we have
$$
|Q_n(-8)\widetilde\Phi^{-1}3^j|
\cdot|\varepsilon_k|
\ge|M''_k|-\biggl|R_n\biggl(-\frac18\biggr)
\widetilde\Phi^{-1}2^{b-3a+1+2j}(-1)^a\biggr|
>\frac13\mspace{2mu},
$$
and so
$$
|\varepsilon_k|>\frac{\widetilde\Phi}{3^{j+1}|Q_n(-8)|}
\ge\frac{\widetilde\Phi}{3^{2\beta}|Q_n(-8)|}\mspace{2mu},
$$
which means that
\begin{equation}
\label{eq16}
|\varepsilon_k|>\frac{\Phi}{3^{2\beta}|Q_n(-8)|}\mspace{2mu}.
\end{equation}
or
\begin{equation}
\label{eq17}
|\varepsilon_k|>\frac{\Phi'}{(n+1)3^{2\beta}|Q_{n+1}(-8)|}\mspace{2mu},
\end{equation}
depending on the choice made in \eqref{eq15}.

\section({A Bound for \$Phi\$}){A Bound for $\Phi$}

For evaluating  $\Phi$ and $\Phi'$ we consider
the functions
\begin{align*}
\varphi(x)&=\min_{0\le y<1}(
-\{-(\alpha+\gamma)x\}
+\{-(\alpha+\gamma)x-y\}+\{y\}
\\
&\qquad\hphantom{=\min_{0\le y<1}(}
-\{(\alpha+\beta+\gamma)x\}
+\{(\alpha+\beta)x+y\}
+\{\gamma x-y\}),
\\
\varphi'(x)&=\min_{0\le y<1}(
-\{(\alpha+\gamma)x+y\}
+\{(\alpha+\gamma)x\}+\{y\}
\\
&\qquad\hphantom{=\min_{0\le y<1}(}
-\{(\alpha+\beta+\gamma)x\}
+\{(\alpha+\beta)x+y\}
+\{\gamma x-y\}),
\end{align*}
which take the values $e_{n,p}$
and $e'_{n,p}$, respectively at the point $m/p$.

All the solutions $x$ of the equation
$\varphi(x)=1$ form the set of intervals
in $[0,1)$, which should contain~$\{x\}$.
If we denote $A_i$ and $B_i$
the left and right points of this intervals, respectively,
then the condition $A_i\le \{m/p\}<\nobreak B_i$ (i.e.\
$e_{n,p}=\varphi(m/p)=1$) is equivalent to
$$
A_i+N\le\frac mp< B_i+N,\qquad N\in\mathbb N,
$$
($\mathbb N$ is the set of non-negative integers),
or the same
$$
p\in\biggl(\frac{m}{B_i+N}\mspace{2mu},
\frac{m}{A_i+N}\biggr],\qquad N\in\mathbb N.
$$
This means that all the prime numbers $p$ such that
$$
p\in\bigcup_{N=0}^{\left[\tfrac{m}{\sqrt{a+b+n}}\right]-1}
\bigcup_i\biggl(\frac{m}{B_i+N}\mspace{2mu},
\frac{m}{A_i+N}\biggr]
$$
(the inequality $p>\sqrt{a+b+n}$ entails $m/(B_i+N)\ge\sqrt{a+b+n}$,
and one can find the bound for $N$)
go to $\Phi$. So we have
\begin{equation}
\label{eq13}
\log\Phi
\ge\sum_{N=0}^{\left[\tfrac{m}{\sqrt{a+b+n}}\right]-1}
\sum_i\biggl(\theta\biggl(\frac{m}{A_i+N}\biggr)
-\theta\biggl(\frac{m}{B_i+N}\biggr)\biggr),
\end{equation}
where $\theta(x)=\sum_{p\le x,\,p\text{ is prime}}\log p$.

The same works for $\varphi'(x)$. And it is proved in \cite{ref6},
that the sets for $\varphi(x)$ and $\varphi'(x)$ differs
only on a set of zero measure.

\section{Analytic and Arithmetical Bounds}
\label{ss5}
Let us take
$$
\alpha = 3,\qquad \beta = 9, \qquad \gamma = 4.
$$

For $\Phi$ we have the set of intervals
$$
\biggl[\frac{1}{8}\mspace{2mu},\frac{1}{7}\biggr]
\cup\biggl[\frac{3}{16}\mspace{2mu},\frac{1}{5}\biggr)
\cup\biggl[\frac{3}{8}\mspace{2mu},\frac{2}{5}\biggr)
\cup\biggl[\frac{9}{16}\mspace{2mu},\frac{4}{7}\biggr]
\cup\biggl[\frac{11}{16}\mspace{2mu},\frac{5}{7}\biggr]
\cup\biggl[\frac{15}{16}\mspace{2mu},1\biggr)
$$
For $\Phi'$ the difference will only be in the right-end
points of the intervals.

We will use the following bounds for
$\theta(x)$ \cite{ref5}: the upper bound
$$
\theta(x)< 1.001102\cdot x\qquad \text{if}\quad x>0,
$$
and the lower bound
$$
\theta(x) > 0.998\cdot x\qquad \text{if}\quad x>487\,381.
$$
Substituting them in \eqref{eq13}
and taking the sum for $N=0,1$ we obtain
\begin{equation}
\label{eqPhi}
\log\Phi>1.639533\cdot m\qquad
\text{for}\quad m\ge 974\,762.
\end{equation}
The same bound holds for $\Phi'$.

We need to estimate the values of $|Q_n(z^{-1})|$, $|R_n(z)|$,
and $|Q_{n+1}(z^{-1})|$, $|R_{n+1}(z)|$ at the point
$z=-1/8$, and to estimate $\Phi'$. We begin with
$$
\log\biggl(\frac{(16m)!}
{(7m-1)!\mspace{2mu}(4m)!\mspace{2mu}
(5m)!}\biggr)\qquad\text{and}\qquad
\log\biggl(\frac{(16m+1)!}
{(7m)!\mspace{2mu}(4m+1)!\mspace{2mu}
(5m-1)!}\biggr).
$$
We use Stirling's formulae \cite{ref3}
$$
\sqrt{2\pi n}\biggl(\frac{n}{e}\biggr)^n
<n!<\sqrt{2\pi n}\biggl(\frac{n}{e}\biggr)^ne^{1/(12n)},
$$
and we find
\begin{align*}
&\log\biggl(\frac{(16m)!}
{(7m-1)!\mspace{2mu}(4m)!\mspace{2mu}
(5m)!}\biggr)
=\log\biggl(\frac{(16m)!}
{(7m)!\mspace{2mu}(4m)!\mspace{2mu}
(5m)!}\biggr)+\log(7m)
\\
&\quad
<\frac12\log(2\pi)+\frac12\log(16)+\frac12\log(m)
+16m\log(16)+16m\log(m)-16m
\\
&\quad\quad+\frac1{12\cdot16m}
\\
&\quad\quad
-\biggl(\frac12\log(2\pi)+\frac12\log(7)+\frac12\log(m)
+7m\log(7)+7m\log(m)-7m\biggr)
\\
&\quad\quad
-\biggl(\frac12\log(2\pi)+\frac12\log(4)+\frac12\log(m)
+4m\log(4)+4m\log(m)-4m\biggr)
\\
&\quad\quad
-\biggl(\frac12\log(2\pi)+\frac12\log(5)+\frac12\log(m)
+5m\log(5)+5m\log(m)-5m\biggr)
\\
&\quad\quad
+\log(7)+\log(m).
\end{align*}
Since
\begin{align*}
&-\log(2\pi)
+\frac12(\log(16)-\log(7)-\log(4)-\log(5))
-\log(m) 
\\
&\qquad\qquad
+\frac1{12\cdot16m}+\log(7)+\log(m)
\\[1mm]
&\qquad <0\qquad\text{for}\quad m\ge974\,762.
\end{align*}
one can have
\begin{align}
\nonumber
&\log\biggl(\frac{(16m)!}
{(7m-1)!\mspace{2mu}(4m)!\mspace{2mu}
(5m)!}\biggr)
\\
&\qquad
<(16\log(16)-7\log(7)-4\log(4)-5\log(5))\cdot m
\label{eqAn}
\\
\nonumber
&\qquad\qquad\qquad\qquad\qquad\qquad\qquad\qquad\qquad\qquad\qquad
\text{for}\quad m\ge974\,762.
\end{align}
Now,
\begin{align*}
&\log\biggl(\frac{(16m+1)!}
{(7m)!\mspace{2mu}(4m+1)!\mspace{2mu}
(5m-1)!}\biggr)
\\
&\qquad
=\log\biggl(\frac{(16m)!}
{(7m)!\mspace{2mu}(4m)!\mspace{2mu}
(5m)!}\biggr)+\log\biggl(\frac{(16m+1)(5m)}{(4m+1)}\biggr)
\\
&\qquad
<\log\biggl(\frac{(16m)!}
{(7m)!\mspace{2mu}(4m)!\mspace{2mu}
(5m)!}\biggr)+\log\biggl(\frac{(16m+4)(5m)}{(4m+1)}\biggr)
\\
&\qquad
=\log\biggl(\frac{(16m)!}
{(7m)!\mspace{2mu}(4m)!\mspace{2mu}
(5m)!}\biggr)+\log(20m),
\end{align*}
and in a similar way we conclude that
\begin{align}
\nonumber
&\log\biggl(\frac{(16m+1)!}
{(7m)!\mspace{2mu}(4m+1)!\mspace{2mu}
(5m-1)!}\biggr)
\\
&\qquad
<(16\log(16)-7\log(7)-4\log(4)-5\log(5))\cdot m +1
\label{eqAn1}
\\
\nonumber
&\qquad\qquad\qquad\qquad\qquad\qquad\qquad\qquad\qquad\qquad\qquad
\text{for}\quad m\ge974\,762.
\end{align}

For the integral in \eqref{eq8p} we write the estimates
\begin{align*}
&\int_0^1t^{4m}(1-t)^{7m-1}(1-zt)^{-16m-1}\,dt
\\
&\qquad
<\biggl(\max_{t\in[0,1]}t^{4}(1-t)^{7}
\biggl(1+\frac t8\biggr)^{-16}\biggr)^{m-1}
\int_0^1t^4(1-t)^{6}\biggl(1+\frac t8\biggr)^{-17}\,dt,
\\
&\int_0^1t^{4m+1}(1-t)^{7m}(1-zt)^{-16m-2}\,dt
\\
&\qquad
<\biggl(\max_{t\in[0,1]}t^{4}(1-t)^{7}
\biggl(1+\frac t8\biggr)^{-16}\biggr)^{m-1}
\int_0^1t^5(1-t)^7\biggl(1+\frac t8\biggr)^{-18}\,dt,
\end{align*}
and so
\begin{align*}
&\log\biggl(\int_0^1t^{4m}(1-t)^{7m-1}(1-zt)^{-16m-1}\,dt\biggr)
\\
&\qquad\qquad\qquad
<-7.884160\cdot (m-1) -8.568400,
\\
&\log\biggl(\int_0^1t^{4m+1}(1-t)^{7m}(1-zt)^{-16m-2}\,dt\biggr)
\\
&\qquad\qquad\qquad
<-7.884160\cdot (m-1) -10.140038.
\end{align*}

Let us check inequality \eqref{eq100} for $n$ and $n+1$:
\begin{align*}
&\log\biggl|R_n\biggl(-\frac18\biggr)
\widetilde\Phi^{-1}2^{b-3a+1+2j}(-1)^a\biggr|
\\
&\qquad
<17.147682\cdot m +\log\biggl(\frac18\biggr)^4\cdot m
-7.884160\cdot (m-1) -8.568400
\\
&\qquad\qquad-1.639533\cdot m+35\log(2)
\\
&\qquad
<-0.693777\cdot m+23.575912
<\log\biggl(\frac23\biggr)\qquad
\text{for}\quad m\ge 974\,762,
\\
&\log\biggl|R_{n+1}\biggl(-\frac18\biggr)
\widetilde\Phi^{-1}2^{b-3a+1+2j}(-1)^a\biggr|
\\
&\qquad<17.147682\cdot m+1+\log\biggl(\frac18\biggr)^4\cdot m
+\log\biggl(\frac18\biggr)-7.884160\cdot (m-1)
\\
&\qquad\qquad
-10.140038-1.639533\cdot m+\log(4)\log(m)+1+35\log(2)
\\
&\qquad
<-0.693777\cdot m +\log(4)\log(m)+21.924832
<\log\biggl(\frac23\biggr)
\\
&\qquad\qquad\qquad\qquad\qquad\qquad\qquad\qquad\qquad\qquad\qquad
\text{for}\quad m\ge 974\,762.
\end{align*}
So, inequality \eqref{eq100} holds,
and we can move on.

For integral in \eqref{eq8q}, in the same
way as for the one in \eqref{eq8p}, one has
\begin{align*}
&\log\biggl(\int_0^1t^{7m-1}(1-t)^{5m}(1+8t)^{4m}\,dt\biggr)
\\
&\qquad\qquad\qquad
<-0.945755\cdot (m-1)-1.725707,
\\
&\log\biggl(\int_0^1t^{7m}(1-t)^{5m-1}(1+8t)^{4m+1}\,dt\biggr)
\\
&\qquad\qquad\qquad
<-0.945755\cdot (m-1)+0.878883.
\end{align*}

Now we can calculate the bounds in \eqref{eq16} and \eqref{eq17}.
We begin with \eqref{eq16}:
\begin{align*}
\log|\varepsilon_k|
&>1.639533\cdot m-18\log(3)-17.147682\cdot m
\\
&\qquad
+0.945755\cdot (m-1)+1.725707
\\
&>-14.562394\cdot m-18.995070
\\
&> -0.81\cdot k
>\log\biggl(\frac49\biggr)\cdot k
\qquad\text{for}\quad k\ge17\,545\,718.
\end{align*}
For \eqref{eq17} we have
\begin{align*}
\log|\varepsilon_k|
&>1.639533\cdot m-\log(4m+1)-18\log(3)-17.147682\cdot m
\\
&\qquad
+0.945755\cdot (m-1)-0.878883
\\
&>-14.562394\cdot m-\log(4)\log(m)-22.599660
\\
&>-14.562414\cdot m -22.599660
\\
&> -0.81\cdot k
>\log\biggl(\frac49\biggr)\cdot k
\qquad\text{for}\quad k\ge17\,545\,718.
\end{align*}

So, we have
$$
\biggl\|\biggl(\frac43\biggr)\biggr\|
>\biggl(\frac49\biggr)^k\qquad\text{for}\quad k\ge17\,545\,718.
$$

\section{The Final check}

We need to check inequality \eqref{eqI} for $6\le k\le 17\,545\,717$.
Following \cite{ref2}, we prove next lemma.

\medskip

\begin{lemma}
\label{lem1}
Let $m$ be a positive integer,
and assume that the number $4^m$
contain no block of $h$ consecutive $0$, or $2$,
in its ternary expansion. Then the inequality
\begin{equation}
\label{main}
\biggl\|\biggl(\frac43\biggr)^k\biggr\|
\ge\biggl(\frac49\biggr)^k
\end{equation}
holds for all
\begin{equation}
\label{main2}
m\biggl(\frac{\log4}{\log9}\biggr)
+\frac{h}2\le k\le m.
\end{equation}
\end{lemma}

\medskip

\begin{proof}
We give a proof by contradiction.
Assume that $k$ is in the specified interval,
but \eqref{main} is not true.
Then for some integer $M_1$ we have one
of the next two equalities:
$$
M_1=\biggl(\frac{4}{3}\biggr)^k\pm\epsilon_1, \qquad
\text{where}\quad 0<\epsilon_1<\biggl(\frac{4}{9}\biggr)^k,
$$
so, with some integer $M_2$
$$
4^m=3^kM_2\mp\epsilon_2, \qquad
\text{where}\quad 0<\epsilon_2<4^m3^{-k}.
$$
Since $m(\log4/\log9)\le k-h/2$, we have
$$
0<\epsilon_2<9^{k-h/2}\cdot3^{-k}=3^{k-h},
$$
but this means, that a block of $h$ digits
of the number $4^m$,
which are responsible for powers
$3^{k-h},3^{k-h+1},\allowbreak \dots,3^{k-1}$,
consists of $2$, or $0$.
\end{proof}

For specified $m$ the software calculates $h(m)$ defined in lemma \ref{lem1},
descends to the new value of $m$ prescribed by \eqref{main2}, and so on.
We started calculations with $m=17\,545\,718$ and stopped at $m=5$.
Results of all the steps are given
in Table~\ref{tab1}.

\begin{table}[ht]
\caption{}\label{tab1}
\renewcommand\arraystretch{1.5}
\noindent\[
\begin{array}{|c|r|r|}
\hline
&\multicolumn{1}{c|}{m}&\multicolumn{1}{c|}{h}
\\ \hline
1 &17545718 &18
\\ \hline 2 &11229269
&16
\\ \hline
3 &7186741
&16
\\ \hline
4 &4599523
&14
\\ \hline
5 &2943702
&14
\\ \hline
6 &1883977
&14
\\ \hline
7 &1205753
&15
\\ \hline
8 &771689
&11
\\ \hline
9 &493886
&13
\\ \hline
10 &316094
&12
\\ \hline
11 &202307
&11
\\ \hline
12 &129482
&10
\\ \hline
13 &82874
&12
\\ \hline
\end{array}
\qquad
\begin{array}{|c|r|r|}
\hline
&\multicolumn{1}{c|}{m}&\multicolumn{1}{c|}{h}
\\ \hline
14 &53046
&11
\\ \hline
15 &33955
&10
\\ \hline
16 &21737
&11
\\ \hline
17 &13917
&8
\\ \hline
18 &8911
&9
\\ \hline
19 &5708
&8
\\ \hline
 20 &3658
&10
\\ \hline
21 &2347
&8
\\ \hline
22 &1507
&7
\\ \hline
23 &968
&5
\\ \hline
24 &622
&7
\\ \hline
25 &402
&10
\\ \hline
26 &263
&5
\\ \hline
\end{array}
\qquad
\begin{array}{|c|r|r|}
\hline
&\multicolumn{1}{c|}{m}&\multicolumn{1}{c|}{h}
\\ \hline
27 &171
&4
\\ \hline
28 &112
&4
\\ \hline
29 &74
&3
\\ \hline
30 &49
&3
\\
\hline 31 &33
&3
\\ \hline
32 &23
&3
\\ \hline
33 &16
&3
\\ \hline
 34 &12
&3
\\ \hline
35 &9
&2
\\ \hline
36 &7
&1
\\
\hline
37 &5
&2
\\
\hline
\multicolumn{3}{c}{}
\\
\multicolumn{3}{c}{}
\\
\end{array}
\]
\end{table}

So Theorem \ref{t1} is proved.

\medskip

The author expresses his gratitude
to Igor P. Rochev for his help in preparing the paper.

\end{document}